\magnification=1200  
\overfullrule=0pt

\baselineskip=14pt

\font\twelvebf=cmbx12

\def\bsh{\backslash}

\def\lgl{\langle}
\def\ot{\otimes}
\def\rgl{\rangle}
\def\sbs{\subset}
\def\sbsq{\subseteq}
\def\sm{\simeq}
\def\wt{\widetilde}

\def\Si{\Sigma}

\def\a{\alpha}
\def\d{\delta}
\def\om{\omega}
\def\vp{\varphi}

\font\tenbb=msbm10
\font\sevenbb=msbm7
\font\fivebb=msbm5
\newfam\bbfam
\textfont\bbfam=\tenbb \scriptfont\bbfam=\sevenbb
\scriptscriptfont\bbfam=\fivebb
\def\bb{\fam\bbfam}

\def\Pb{{\bb P}}

\def\build#1_#2^#3{\mathrel{
\mathop{\kern 0pt#1}\limits_{#2}^{#3}}}

\def\ra{\rightarrow}
\def\hra{\hookrightarrow}

\def\Im{\mathop{\rm Im}\nolimits}
\def\Ker{\mathop{\rm Ker}\nolimits}

\def\rk{\mathop{\rm rk}\nolimits}

\def\Oc{{\cal O}}

\def\XX{\hbox{\vrule\vbox to8pt{\hrule
width 4pt \vfill\hrule}\vrule}}


\centerline{\twelvebf Appendix: On linear subspaces contained 
in}

\centerline{\twelvebf the secant varieties of a projective 
curve}

\medskip

\centerline{\bf Claire Voisin}

\smallskip

\centerline{\sl Institut de Math\'ematiques de Jussieu}

\centerline{\sl CNRS, UMR 7586}

\vglue 1cm

\noindent {\bf 1. Introduction.} 

\smallskip

If $C \sbs \Pb^N$ is a curve 
imbedded in projective space, one can consider the secant variety 
$\Si_d = \build\cup_{Z \, \in \, C^{(d)}}^{} \lgl Z \rgl$ swept out 
by the linear spans of $d$-uples of points of $C$. This $\Si_d$ 
contains the $\Pb^{d-1}$'s parametrized by $Z \in C^{(d)}$ (here we 
are assuming that $d$ is not large with respect to $N$).
More precisely, $\Si_d$ is birational to a projective bundle of rank
$d-1$ over $C^{(d)}$. On the other hand,
 if $d$ is large enough, $C^{(d)}$ also
contains
positive dimensional
projective spaces, corresponding to linear systems on $C$. Deciding
whether or 
not $\Si_d$ contains  linear subspaces other than those contained
in some of the $\Pb^{d-1}_Z$'s is thus a non trivial problem.

Some time ago, C. Soul\'e obtained estimates for the maximal dimension of a 
linear subspace contained in $\Si_d$, and asked me whether
an ad hoc geometric argument would lead to  other results.

  One answer in this direction is as follows:

\smallskip

We assume that $C$ is smooth of genus $g > 0$ and that the embedding 
$C \sbs \Pb^N$ is given by the sections of a line bundle $L \ot 
\om_C$, with ${\rm deg}(L) = m$. We then show:

\medskip

\noindent {\bf Theorem.} {\it If $m \geq 2d + 3$, and $\d \geq d-1$, 
any $\Pb^{\d}$ contained in $\Si_d$ is one of the $\Pb^{d-1} = \lgl Z 
\rgl$, $Z \in C^{(d)}$. In particular, $\Si_d$ contains no projective 
space $\Pb^{\d}$, for $\d \geq d$.}

\bigskip

\noindent {\bf Thanks.} I wish to thank Christophe Soul\'e for 
interesting discussions and for providing the motivation to write 
this Note.

\vglue 1cm

\noindent {\bf 2. Proof of the theorem.} 

\smallskip

We first recall a few basic facts about secant varieties of curves 
(see [1]). First of all, since $m \geq 2d + 1$, for any 
effective divisor $Z$ of degree $k \leq 2d$ on $C$, we have $H^1 (L 
\ot \om_C (-Z)) = 0$, hence the linear span of $Z$ is of dimension 
$k-1$. Let now $E \ra C^{(d)}$ be the vector bundle with fiber $H^0 
(L \ot \om_{C \mid Z})$ at $Z \in C^{(d)}$. Since the restriction map 
$H^0 (L \ot \om_C) \ra H^0 (L \ot \om_{C \mid Z})$ is surjective for 
any $Z \in C^{(d)}$, there is a well defined morphism $\a : \Pb 
(E^*) \ra \Pb^N$, whose image is exactly the secant variety $\Si_d$. 
Since sections of $L \ot \om_C$ separates any $2d$ points on $C$, it 
follows that $\a$ is one to one over $\Si_d - \Si_{d-1}$. An easy 
computation shows that for any $Z \in C^{(d)}$, and for any $x$ in 
the linear span of $Z$, but not in the linear span of any $Z' 
\not\sbsq Z$, the differential of $\a$ is of maximal rank, so that 
$\Si_d \bsh \Si_{d-1}$ is smooth of dimension $2d-1$. The 
projectivized tangent space to $\Si_d$ at $\a (x)$ is easy to 
describe, at least when $Z$ is a reduced divisor $\build\sum_{1}^{d} 
z_i$: indeed this is a $\Pb^{2d-1}$ which contains $\lgl Z \rgl$ and 
also each projective line tangent to $C$ at some point $z_i \in Z$, 
as one sees by deforming $Z$ fixing $z_j$, $j \ne i$. It follows that 
it must be equal to the linear span of the divisor $2Z$. By 
continuity, this description of the projectivized tangent space to 
$\Si_d$ remains true at any point of $\Si_d - \Si_{d-1}$ .

\smallskip

We now start the proof of the theorem. We suppose that $\d \geq d-1$, 
and assume that some projective space $\Pb^{\d}$ is contained in 
$\Si_d$. Assuming $\Pb^{\d}$ is not contained in one of the 
$\Pb^{d-1}_Z$'s we shall derive a contradiction.

\smallskip

Note that by induction on $d$, we may assume that $\Pb^{\d}$ is not 
contained in $\Si_{d-1}$. Let $\wt{\Pb}^{\d}$ be the closure of 
$\a^{-1} (\Pb^{\d} \bsh \Pb^{\d} \cap \Si_{d-1})$ in $\Pb (E^*)$. 
Denote by $\pi : \wt{\Pb}^{\d} \ra C^{(d)}$ the restriction to 
$\wt{\Pb}^{\d}$ of the structural projection $\Pb (E^*) \ra C^{(d)}$. 
Let $W := \pi \, (\wt{\Pb}^{\d})$ and $w := \dim W$. Our assumption 
is that $w > 0$. We shall denote by $P_v$ the fiber $\pi^{-1} (v)$. 
It is a projective space $\Pb^{\d} \cap \lgl Z_v \rgl$, which is 
generically of dimension $s = \d - w$.

\smallskip

We start with the following observation:

\medskip

\noindent {\bf Lemma 1.} {\it Under our assumption $\dim W > 0$ we 
have the inequality
$$
w > \d - w \, . \leqno (1)
$$
}

\medskip

\noindent {\it Proof.} Indeed, we may assume that for $v,v'$ two 
generic distinct points of $W$, the supports of the associated 
divisors $Z_v$, $Z_{v'}$ of $C$ are disjoint. Otherwise, $Z_v$ would 
contain a fixed point $x \in C$, for any $v \in W$. But projecting 
$C$ from $x$, we then get a curve $C' \sbs \Pb^{N-1}$, such that 
$\Si'_{d-1}$ contains a $\Pb^{\d - 1}$ which is not a $\Pb_Z^{d-2}$; 
since we may assume the theorem proven for $(m-1 , d-1)$, this is 
impossible.

\smallskip

Now choose $v,v'$ as above. The projective spaces $\lgl Z_v \rgl$ and 
$\lgl Z_{v'} \rgl$ do not meet, hence the projective spaces $P_v = 
\lgl Z_v \rgl \cap \Pb^{\d}$, $P_{v'} = \lgl Z_{v'} \rgl \cap 
\Pb^{\d}$ do not meet. Since they are of dimension $s$ in a 
$\Pb^{\d}$, it follows that $2s < \d$, or $w > \d - w$. \hfill \XX

\bigskip

Next we observe that, at each point $\a (x,Z)$ of $\Pb^{\d} - 
(\Pb^{\d} \cap \Si_{d-1})$, $\Pb^{\d}$ is contained in the 
projectivized tangent space of $\Si_d$ at $\a (x,Z)$, that is in 
$\lgl 2Z \rgl$. Hence for any $v \in W$, the corresponding divisor 
$Z_v \in C^{(d)}$ satisfies
$$
\Pb^{\d} \sbs \lgl 2 \, Z_v \rgl \, .
$$
We next study the infinitesimal variation of $\lgl 2 \, Z_v \rgl \sbs 
\Pb^N$. Let $H := \Oc_{\Pb^N} (1)$. Then we have the identification
$$
H^0 (\Pb^N , H) \sm H^0 (C , L \ot \om_C) \, , \leqno (2)
$$
which by definition of the linear span, induces an identification
$$
H^0 (\Pb^N , H \ot I_{\lgl 2 Z_v \rgl}) \sm H^0 ( C , L \ot \om_C (-2 
\, Z_v)) \, . \leqno (3)
$$
If $h \in T_{W,v}$, the infinitesimal deformation of $\lgl 2 \, Z_v 
\rgl$ in the direction $h$ is described by an homomorphism:
$$
\vp_h : H^0 (\Pb^N , H \ot I_{\lgl 2 Z_v \rgl}) \ra H^0 (\lgl 2 \, 
Z_v \rgl , H_{\mid \lgl 2 Z_v \rgl}) \, .
$$
We have now an isomorphism induced by (2) and (3):
$$
H^0 (\lgl 2 \, Z_v \rgl , H_{\mid \lgl 2 Z_v \rgl}) \sm H^0 (L \ot 
\om_{C \mid 2 Z_v}) \, . \leqno (4)
$$
We have the following

\medskip

\noindent {\bf Lemma 2.} {\it Under the isomorphisms $(3)$ and $(4)$, 
if we identify $h$ to an element $u_h \in H^0 (\Oc_C (Z_v)_{\mid 
Z_v})$, $\vp_h$ identifies to the multiplication
$$
u_h : H^0 (C , L \ot \om_C (-2 \, Z_v)) \ra H^0 (Z_v , L \ot \om_C 
(-Z_v)_{\mid Z_v})
$$
followed by the inclusion
$$
H^0 (Z_v , L \ot \om_C (-Z_v)_{\mid Z_v}) \hra H^0 (2 \, Z_v , L \ot 
\om_{C \mid 2 Z_v}) \, . 
$$
}

\medskip

The proof is straightforward once we recall the construction of 
$\vp_h$ by differentiating under the parameters the equations 
vanishing on $\lgl 2 \, Z_v \rgl$. \hfill \XX

\bigskip

We know that the spaces $\lgl 2 \, Z_v \rgl$, for $v \in W$, contain 
$\Pb^{\d}$. Infinitesimally, this translates into the fact that for 
any $h \in T_{W,v}$, the image of $\vp_h$ vanishes on $\Pb^{\d}$, 
that is, is contained in
$$
\Ker (H^0 (\lgl 2 \, Z_v \rgl , H_{\mid \lgl 2 Z_v \rgl}) \ra H^0 
(\Pb^{\d} , H_{\mid \Pb^{\d}})) \, .
$$

>From the description of $\vp_h$ given in Lemma~2, we see that $\Im \, 
\vp_h$ is contained in
$$
K := \Ker (H^0 (\lgl 2 \, Z_v \rgl , H_{\mid \lgl 2 Z_v \rgl}) \ra 
H^0 (\lgl Z_v \rgl , H_{\mid \lgl Z_v \rgl})) \, .
$$
Indeed, via the isomorphism (4), $K$ identifies to
$$
\Ker (H^0 (L \ot \om_{C \mid 2 Z_v}) \ra H^0 (L \ot \om_{C \mid 
Z_v})) = \Im \, H^0 (L \ot \om_C (-Z_v)_{\mid Z_v}) \ra H^0 (L \ot 
\om_{C \mid 2 Z_v}) \, .
$$
Finally, note that the restriction map $K \ra H^0 (\Pb^{\d} , 
H_{\mid \Pb^{\d}})$ has rank equal to the dimension of
$$
\Ker (H^0 (\Pb^{\d} , H_{\mid \Pb^{\d}}) \ra H^0 (\Pb^{\d} \cap \lgl 
Z_v \rgl , H_{\mid \Pb^{\d} \cap \lgl Z_v \rgl}))\, ,
$$
which is equal to $\d - s$, since $\Pb^{\d} \cap \lgl Z_v \rgl = P_v$ 
is of dimension $s$.

 Denote now by $V \sbs H^0 (\Oc_C (Z_v)_{\mid Z_v})$ 
the tangent space to $W$ at $v$. Lemma~2 and the estimate above give 
us the following conclusion:

\medskip

\noindent {\bf Lemma 3.} {\it Under our assumptions, the 
multiplication map
$$
\mu : V \ot H^0 (C , L \ot \om_C (-2 \, Z_v)) \ra H^0 (L \ot \om_C 
(-Z_v)_{\mid Z_v})
$$
has its image contained in a subspace of codimension at least $w$.} 
\hfill \XX

\bigskip

We now derive a contradiction. We observe first that since 
$\wt{\Pb}^{\d}$ is a rational variety dominating $W$, $W$ is 
contained in a linear system $\vert D \vert \sbs C^{(d)}$. Hence 
$\Oc_C (D) = \Oc_C (Z_v)$ for all $ v \in W$, and the fact that $W \sbs 
\vert D \vert$ translates infinitesimally into the fact 
that $V = T_{W,v}$ is 
contained in the image of the restriction map:
$$
H^0 (\Oc_C (Z_v)) \ra H^0 (\Oc_C (Z_v)_{\mid Z_v}) \, .
$$
Let now $\wt V$ be the inverse image of $V$ under this restriction 
map. Then $\rk \, \wt V = w+1$, and Lemma~3 shows that the 
multiplication map
$$
\wt\mu : \wt V \ot H^0 (C , L \ot \om_C (-2 \, Z_v)) \ra H^0 (C , L 
\ot \om_C (-Z_v))
$$
has its image contained in a space of codimension at least $w$.

\smallskip

Now we have the equality:
$$
\rk \, H^0 (C , L \ot \om_C (-Z_v)) = d + \rk \, H^0 (C , L \ot \om_C 
(-2 \, Z_v)) \, ,
$$
since $H^1 (C , L \ot \om_C (-2 \, Z_v)) = 0$. So we conclude that
$$
\rk \, \wt\mu \leq h^0 (C , L \ot \om_C (-2 \, Z_v)) + d - w \, . 
\leqno (5)
$$
On the other hand, we can apply Hopf lemma to $\wt\mu$, and the 
inequality in Hopf lemma must be strict here, since the line bundle 
$L \ot \om_C (- 2 \, Z_v)$ is very ample, being of degree at least 
$2g + 1$, and $C$ is not rational. This gives us:
$$
\rk \, \wt\mu > w + 1 + h^0 (C , L \ot \om_C (-2 \, Z_v)) - 1 \, .
\leqno (6)
$$
Combining (5) and (6), we get:
$$
d-w > w \, . \leqno (7)
$$
But this contradicts inequality (1), since $\d \geq d-1$. \hfill \XX

{\bf References.} 

[1]  A. Bertram : Moduli of rank 2 vector bundles,
theta divisors, and the geometry of curves in projective space, {\it
J. Diff. Geom}. 35, 1992, 429-469.

\bye